%% file: main.tex
\documentclass[1p]{elsarticle}
\newcommand{\dumblang}[3]{#2} 
\usepackage{mathtools}
\mathtoolsset{showonlyrefs}
\usepackage{hyperref}
\usepackage{amsmath, amssymb, amsfonts, amsthm}
\usepackage{mathtools,stmaryrd,eucal,marvosym}
\usepackage{bm}
\makeatletter
\def\ps@pprintTitle{%
  \let\@oddhead\@empty
  \let\@evenhead\@empty
  \let\@oddfoot\@empty
  \let\@evenfoot\@oddfoot
}
\makeatother
\newtheoremstyle{mystyle}
  {}
  {}
  {\itshape}
  {}
  {\bfseries}
  {.}
  { }
  {\thmname{#1}\thmnumber{ #2}\thmnote{ (#3)}}
\theoremstyle{mystyle}
\newtheorem{satz}{Theorem}[section]
\newtheorem{lem}[satz]{Lemma}

\newenvironment{bew}[1][\proofname]{%
  \proof[\itshape\bfseries #1\hspace{-0.36em} ]%
}{\endproof}
\theoremstyle{definition}

\newtheorem{defi}[satz]{Definition}
\newtheorem{rem}[satz]{Remark}

\newtheorem{beisp}[satz]{\dumblang{Beispiel}{Example}{Voorbeeld}}

\usepackage[english]{babel}
\begin{document}
\title{
Stochastic Delay Differential Equations have blow-up solutions if and only if their
instantaneous counterparts have them
}
\date{\today}
\author[1,2]{Julius Busse}
\ead{julius.busse@ucdconnect.ie}
\affiliation[1]{organization={Mathematical Institute, Leiden University},
    city={Leiden},
    country={The Netherlands}}
\affiliation[2]{organization={School of Mathematics and Statistics, University College Dublin},
    city={Belfield, Dublin 4},
    country={Ireland}}

\begin{abstract}
\noindent
Motivated by a recent publication by Ishiwata and Nakata \cite{jap}, we
prove that sufficiently regular stochastic delay differential equations (SDDEs) with
a single discrete delay have blow up solutions
if and only if their undelayed counterparts have them, using a comparison theorem by
Ikeda and Watanabe \cite{watanabe}. This result has applications
in mathematical biology and finance.
\end{abstract}
\maketitle
\input{00_intro.tex}
\input{20_math_basics.tex}
\input{40_delay.tex}
\input{60_summary.tex}
\bibliographystyle{alpha}
\bibliography{mybibliography}
\end{document}

%% file: 00_intro.tex
\section{Introduction}
Delay differential equations appear in many applications, such as in population
modeling \cite{rihan}, modeling the glucose-insulin system or other processes in the
human body \cite{glucose-insulin} and in control problems in engineering
\cite{delay_engineering}. We introduce stochasticity to these models.
``All biological dynamical systems
evolve under stochastic forces, if we define stochasticity as the parts of the dynamics that we
either cannot predict or understand or that we choose not to include in the explicit
modeling'' \cite[p.3]{biomath}. The resulting SDDEs are a growing research interest,
both from a purely mathematical perspective \cite{preprint_onno}, in mathematical
biology \cite{ditlevsen}, and also in mathematical finance, where they
can be used to model insider trading \cite{leeDelayFinancialModel2011}.

Explosion properties of these SDDEs are of interest because of their applications in
mathematical biology . In mathematical biology, there is a research
interest in SDE models with non locally Lipschitz noise coefficients
\cite{birmingham}, \cite[p.12]{biomath}. Using the It\^o formula, we can transform
the models we study onto \( (0,\infty) \), where explosion becomes extinction (see Example
\ref{ex:delay}) and the noise coefficients are non locally Lipschitz, thus allowing us to gain insight about explosion properties and
estimates of these probabilities using our results (more in \cite{masters}). Non
locally Lipschitz noise coefficients are also used in mathematical finance, for
example when modelling how interest rates evolve over time
\cite{coxTheoryTermStructure}, where a realistic restriction is that interest rates
may never be negative. Again this can be answered by transforming models we study
onto \( (0,\infty)\) using the It\^o formula.

We combine the comparison theorem by Ikeda and Watanabe as presented in \cite{watanabe} and the
analysis of blow-up solutions for scalar differential equations with discrete delay
as described in \cite{jap} to reduce explosion behavior of sufficiently regular SDDE
models to that of SDEs. Through an application of It\^o's formula, this can then be
applied to extinction properties of delayed population models with non locally Lipschitz
noise coefficients \cite{masters}.

%% file: 20_math_basics.tex
\section{Mathematical Basics}
\subsection{Basics of SDEs}
\subsubsection{Elementary definitions (SDEs)}
For basic notions from probability theory, including Brownian motion,
we refer to \cite{karatzas_shreve}. Analogously to the definition of a strong solution given earlier in \cite[p.285]{karatzas_shreve}, we can also define the
concept of a strong solution up to an explosion time. Thus, we arrive
at the following definition.

\begin{defi}[Strong local solutions of SDEs]
	A strong solution up to an explosion time of the SDE \begin{equation}
		\,\mathrm{d}x(t)=b(x(t),t)\,\mathrm{d}t+\sigma(x(t),t)\,\mathrm{d}W(t)
	\end{equation}
	on the given probability space \((\Omega,\mathcal{F},\mathbb{P})\)
	with filtration \((\mathcal{F}_t)_{t\ge 0}\) and
	with respect to the fixed Brownian motion \(W\) on
	\((\Omega,\mathcal{F},\mathbb{P}),(\mathcal{F}_t)_{t\ge 0}  \) and initial condition
	\(x_0\) (a \(\mathcal{F}_0\)-measurable random variable) is a process \(x=\{x(t) : 0\le t<\infty{}\} \) with
	sample paths that are continuous in the extended real numbers
	\(\mathbb{R}\cup\{-\infty{},+\infty{}\}\) and with the following properties.
	\begin{enumerate}[(i)]
		\item
		      \(x\) is adapted to the filtration \((\mathcal{F}_t)_{t\ge 0}\),
		      \(\mathbb{P}[x(0)=x_0]=1\).
		\item
		      Define \(S_n:=\inf\{t\ge 0:|x(t)|>n\}\) (possibly \(\infty{}\)),
		      then we have \begin{equation}
			      \mathbb{P}\left[\int_{0}^{t\land{}S_n}|b(x(s),s)|+\sigma^2(x(s),s)\,\mathrm{d}s<\infty{}   \right]
			      =1
		      \end{equation}
		      for all \(0\le t<\infty{}\), and we have
		\item
		      \begin{align}
                  x(t\land{}S_n)=&x(0)+\int_{0}^{t}b(x(s),s)1_{\{s\le
                  S_n\}}\,\mathrm{d}s\\
                                 &+\int_{0}^{t}\sigma(x(s),s)1_{\{s\le S_n\}}\,\mathrm{d}W(s)
		      \end{align}
		      for all \(t\ge 0\), \(n\ge 1\) almost surely. We refer to
		      \begin{equation}
			      S=\lim_{n\to \infty{}}S_n
		      \end{equation}
		      as the \emph{explosion time} for \(x\). We say \(x\) explodes or \(x\)
		      blows up at time \(S\).
	\end{enumerate}
\end{defi}

\begin{rem}
	If explosion occurs towards \(+\infty{}\), we denote the corresponding
	explosion time by \(S^+\)
	or \(T_{+\infty}\), and if it occurs towards \(-\infty{}\), we denote
	it by \(S^-\) or \(T_{-\infty}\).
	Similarly, in the case of \(x_0>0\), we call the time
	\((S_0=)T_0:=\lim_{n\to \infty{}}\inf\{t\ge 0:x(t)\notin
	(\frac{1}{n},\infty{})\}\) the extinction time of \(x\).
\end{rem}

\subsubsection{Existence and uniqueness results (SDEs)}
The following definition is taken from \cite[p.250]{protter}
\begin{defi}[(Locally) Lipschitz functions]
	A function \(f:\mathbb{R}\cup\{-\infty{},\infty{}\}\times \mathbb{R}_+\to \mathbb{R}\) is
	\emph{Lipschitz} if there exists a constant \(k\) such that
	\begin{itemize}
		\item
		      \(|f(x,t)-f(y,t)|\le k|x-y|\), for each \(t>0\),
		      \(x,y\in\mathbb{R}\cup\{-\infty{},\infty{}\}\), and
		\item
		      \(t\mapsto f(t,x)\) is right continuous with left limits, for each
		      \(x\in \mathbb{R}\cup\{-\infty{},\infty{}\}\).
	\end{itemize}
	Similarly, \(f\) is called \emph{locally Lipschitz} if for every
	compact \(K\subset{}\mathbb{R}\cup\{-\infty{},\infty{}\}\) there exists a \(k\) such that
	\(|f(x,t)-f(y,t)|\le k|x-y|\) for all \(x,y\in K\) and all \(t>0\).
\end{defi}

The following result is mentioned in \cite[p.254]{protter}.
\begin{satz}[Existence and uniqueness of SDE solutions]
	\label{thm:ex_and_uniq}
	Given an SDE \begin{equation}
		\,\mathrm{d}x(t)=b(x(t),t)\,\mathrm{d}t+\sigma(x(t),t)\,\mathrm{d}W(t),
	\end{equation}
	where \(b\) and \(\sigma \) are locally Lipschitz as defined above.
	Then there exists a unique strong solution up to an explosion time \(T\).
\end{satz}

\subsection{Basics of SDDEs}
\subsubsection{Elementary definitions (SDDEs)}
For the definition of a solution of an SDDE and for the corresponding
existence and uniqueness results we refer to
\cite{preprint_onno}.
Fix a probability space \((\Omega,\mathcal{F},\mathbb{P})\) together
with a filtration \((\mathcal{F}_t)_{0\le t<\infty{}}\) satisfying the usual conditions.
Denote by \(\mathbb{D}[0,\infty{})\) the space of
\(\mathcal{F}_t\)-adapted c\`adl\`ag processes \(\{x(t):0\le
t<\infty{}\}\). One can trivially extend the filtration by setting
\(\mathcal{F}_s=\mathcal{F}_0\) for \(-\tau\le s<0\), for some
\(\tau>0\) fixed, and in a similar fashion to the above we write
\(\mathbb{D}[-\tau,\infty{})\) for the space of
\((\mathcal{F}_t)_{-\tau\le t<\infty{}}\)-adapted c\`adl\`ag processes
\(\{x(t):-\tau\le t<\infty{}\}\). The space \(\mathbb{D}[-\tau,0]\)
consists of all c\`adl\`ag processes \(\{x(t):-\tau\le t\le 0\}\)
where \(x(s)\) is \(\mathcal{F}_0\)-measurable for all times
\(s\in[-\tau,0]\).

A stochastic process \(x=\{x(t):-\tau\le
t<\infty{}\in\mathbb{D}[-\tau,\infty{})\}\) gives rise to the segment
process \(\{x_t:t\ge 0\}\), where each segment \(x_t\) is a
\(\mathcal{F}_t\)-measurable c\`adl\`ag process on \([-\tau,0]\)
defined by \begin{equation}
	x_t(\theta):=x(t+\theta),\quad -\tau\le \theta\le 0.
\end{equation}
The segment process can be regarded as an \(\mathcal{F}_t   \)-adapted
\(\mathbb{D}[-\tau,0]   \)-valued stochastic process, thus where each
\(x_t\) is a \(\mathbb{D}[-\tau,0]\)-valued random variable.
Let \(W\) denote a Brownian motion. We wish to provide a rigorous
mathematical interpretation for autonomous equations of the form \begin{equation}
	\begin{cases}
		\,\mathrm{d}x(t)=a(x_t)\,\mathrm{d}t+b(x_{t-})\,\mathrm{d}W(t),\quad & \text{for
		}t\ge 0,                                                                                         \\
		x(u)=\Phi(u),\quad                                                   & \text{for }u\in[-\tau,0],
		\label{eq:SDDE_definition}
	\end{cases}
\end{equation}
subject to an initial process \(\Phi\in\mathbb{D}[-\tau,0]\).

Equation \eqref{eq:SDDE_definition} is shorthand notation for the
integral equation \begin{equation}
	x(t)=\Phi(0)+\int_{0}^{t}\bm A_\Phi(x)(s)\,\mathrm{d}s+\int_{0}^{t}\bm
	B_\Phi(x)(s-)\,\mathrm{d}W(s),
	\label{eq:SDDE_def_integral_formulation}
\end{equation}
\( t\ge 0 \),
where
\begin{enumerate}[(i)]
	\item
	      the maps \(\bm A_\Psi,\bm B_\Psi:\mathbb{D}[0,\infty{})\to
	      \mathbb{D}[0,\infty{})\), for some fixed \(\Psi\in\mathbb{D}\), are
	      defined pathwise for all \(y=\{y(t):t\ge
	      0\}\in\mathbb{D}[0,\infty{})\) by \begin{align}
              &\bm A_\Psi(y)(s,\omega)=a(y_s^\Psi(\omega))\text{  and  }\\
              &\bm B_\Psi(y)(s,\omega)=b(y_s^\Psi(\omega)),\ s\ge 0;
	      \end{align}
	\item
	      the \((y_t^\Psi)_{t\ge 0}\) is the segment process of
	      \(y^\Psi\in\mathbb{D}[-\tau,\infty{})\), where \begin{equation}
		      y^\Psi(s)
		      =
		      \begin{cases}
			      \Psi(s), & -\tau\le s<0, \\
			      y(s),    & s\ge 0;
		      \end{cases}
	      \end{equation}
	\item
	      the functional \(a:\mathbb{D}[-\tau,0]\to\mathbb{R}\) and
	      \(b:\mathbb{D}[-\tau,0]\to\mathbb{R}\) will be assumed to be locally
	      Lipschitz and continuous.
\end{enumerate}

\begin{defi}[Strong local solution of SDDE]
	\label{def:sdde}
	Let \(x=\{x(t):-\tau\le t<\infty{}   \}\) be a stochastic process and
	\(S_\infty{}\) be a stopping time. Then \((x,S_\infty{})\), often
	abbreviated by \(x\) again, is said to be a strong solution up to an
	explosion time of the equation \eqref{eq:SDDE_definition}
	on the interval \([-\tau,S_\infty{})\), if \(x_0=\Phi \) holds
	\(\mathbb{P}\)-a.s., \(x^{S_k}\in\mathbb{D}[0,\infty{})\), and we have
	\(\mathbb{P}\)-a.s. \begin{align}
		x(t\land{}T_k)=\Phi(0)
        &\int_{0}^{t\land{}S_n}\bm A_\Phi(x^{S_n})(s)\,\mathrm{d}s+\\
        +&\int_{0}^{t\land{}S_n}
		\bm B_\Phi(X^{S_n})(s-)\,\mathrm{d}W(s),
	\end{align}
    \( t\ge{0} \),
	where \((S_n)_{n\ge 1}\) is a non-decreasing sequence of finite
	stopping times such that \(S_n\uparrow S_\infty{}\) holds
	\(\mathbb{P}\)-a.s. as \(n\to \infty{}\).
\end{defi}

\subsubsection{Existence and uniqueness results (SDDEs)}
\begin{defi}[Locally Lipschitz functional]
	A functional \(f:C[-\tau,0]\to\mathbb{R}\) is called \emph{locally
		Lipschitz} if for every \(r>0\) there exists an \(L=L(r)>0\), such
	that for all \(x,y\in C[-\tau,0]\) for which
	\(\|x\|_\infty{},\|y\|_\infty{}<r\) we have \begin{equation}
		|f(x)-f(y)|<L\|x-y\|_\infty{}.
	\end{equation}
\end{defi}

\begin{satz}[Existence and uniqueness of SDDE solutions]
	\label{thm:SDDE_exist}
	Let \(a,b:C[-\tau,0]\to\mathbb{R}\) be locally Lipschitz functionals.
	Consider the SDDE \begin{equation}
		\,\mathrm{d}x(t)=a(x_t)\,\mathrm{d}t+b(x_t)\,\mathrm{d}W(t).
		\tag{\eqref{eq:SDDE_definition}}
	\end{equation}
	Then for every \(\mathcal{F}_0\)-measurable random variable \(\Phi \)
	with values in \(C[-\tau,0]\), there exists a stopping time
	\(S_\infty{}\) and process \(x\) with continuous sample paths on
	\([-\tau,S_\infty{})\), such that \((x,S_\infty{})\) is a local
	solution to equation \eqref{eq:SDDE_definition} with \(x_0=\Phi \).
	Moreover, the stopping time \(S_\infty{}\) can be chosen such that
	\((x,S_\infty{})\) is a maximal solution, which is also unique.
\end{satz}

\begin{rem}
	We often consider equations with discrete delay, i.e.\ equations of the
	form
	\begin{equation}
		\,\mathrm{d}x(t)=[b(x(t))+c(x(t-\tau))]\,\mathrm{d}t+\sigma(x(t))\,\mathrm{d}W(t).
	\end{equation}
	This should be interpreted as \begin{align}
		\,\mathrm{d}x(t) & =a(x_t)\,\mathrm{d}t+\bm{\sigma}(x_t)\,\mathrm{d}W(t),
		\intertext{where}
		a(x_t)           & =b(x_t(0))+c(x_t(-\tau)),                              \\
		\bm{\sigma}(x_t) & =\sigma(x_t(0)).
	\end{align}
\end{rem}

\subsection{A comparison theorem for SDEs}

The following theorem is taken from \cite[Theorem 1.1]{watanabe}. It
is of relevance to sharper results on the explosion behavior of SDDEs
we prove later on, because it also covers the nonautonomous case.

\begin{satz}[Comparison theorem by Ikeda and Watanabe]
	\label{thm:watanabe_comparison_sde}
	Suppose we are given the following:
	\begin{enumerate}[(i)]
		\item
		      A real continuous function \(\sigma\) defined on
		      \([0,\infty{})\times \mathbb{R}\) such that \begin{equation}
			      |\sigma(x,t)-\sigma(y,t)|\le \rho(|x-y|),\quad x,y\in\mathbb{R},\
			      t\ge0,
			      \label{eq:K}
		      \end{equation}
		      where \(\rho\) is an increasing function on \([0,\infty{})\) such
		      that \(\rho(0)=0\) and \begin{equation}
			      \lim_{\epsilon\downarrow 0}\int_{0}^\epsilon\rho(\xi)^{-2}\,\mathrm{d}\xi=\infty{},
		      \end{equation}
		\item
		      real continuous functions \(b_1\) and \(b_2\) defined on
		      \(\mathbb{R}\times[0,\infty{})\) such that \begin{equation}
			      b_1(x,t)<b_2(x,t),\quad t\ge 0,\ x\in \mathbb{R}.
			      \label{watanabe:1.8}
		      \end{equation}
		      Let \((\Omega,\mathcal{F},\mathbb{P})\) be a complete
		      probability space, \((\mathcal{F}_t)_{t\ge 0}\) a filtration
		      satisfying the usual conditions, and suppose we are
		      given the following stochastic processes defined on it:
		\item
		      two real \((\mathcal{F}_t)\)-adapted continuous processes
		      \((x_1(t))_{t\ge 0}\) and \((x_2(t))_{t\ge 0}\),
		\item
		      a one-dimensional \((\mathcal{F}_t)\)-Brownian motion \(B\)
		      such that \(B(0)=0\) a.s.,
		\item
		      two real \(\mathcal{F}_t\)-adapted measurable processes
		      \((\beta_1(t))_{t\ge 0}\) and \((\beta_2(t))_{t\ge 0}\).
	\end{enumerate}
	We assume that they satisfy the following conditions with probability
	one:
	\begin{align}
		 &
         x_i(t)-x_i(0)=\int_{0}^{t}\sigma(x_i(s),s)\,\mathrm{d}B(s)+\int_{0}^{t}\beta_i(s)\,\mathrm{d}s,\\
         &i=1,2\text{ for all }t\ge 0,                                                                            \\
		 & x_1(0)\le x_2(0),                                                                                    \\
		 & \beta_1(t)\le b_1(x_1(t),t),\text{ for all }t\ge 0,                                                  \\
		 & \beta_2(t)\ge b_2(x_1(t),t),\text{ for all }t\ge 0,
	\end{align}
	then, with probability one, we have \begin{equation}
		x_1(t)\le x_2(t),\text{ for all }t\ge 0.
	\end{equation}
	If furthermore at least of the following SDEs
	\begin{equation}
		\,\mathrm{d}x(t)=\sigma(x(t),t)\,\mathrm{d}B(t)+b_i(x(t),t)\,\mathrm{d}t,\quad
		i=1,2,
	\end{equation}
	has a unique solution,
	then we have the same conclusion by weakening
	\eqref{watanabe:1.8} to \begin{equation}
		b_1(x,t)\le b_2(x,t),\quad t\ge 0,\ x\in \mathbb{R}.
	\end{equation}
\end{satz}
\subsection{Coupling theory}
The following two definitions are taken from den Hollander's
lecture notes on coupling theory, \cite[p.10f]{coupling}. Coupling
theory is a technique used to study and compare the behavior of
probability measures and random variables by constructing them on a
common probability space. The name-giving idea is to define a
\textbf{coupling}, which is a probability measure/pair of random
variables on the product space whose marginals are identical to/have
the same distribution as the original objects. Coupling is an
important tool in Markov chain analysis, the study of random walks and
interacting particle systems. Den Hollander provides an overview of
these cases in \cite{coupling}.

\begin{defi}[Coupling of probability measures]
	A coupling of two probability measures \(\mathbb{P}\) and
	\(\mathbb{P}'\) on the same measurable space \((\Omega,\mathcal{F})\)
	is any probability measure \(\widehat{\mathbb{P}}\) on the product
	measurable space
	\((\Omega\times\Omega,\mathcal{F}\otimes{}\mathcal{F})\) (where
	\(\mathcal{F}\otimes{}\mathcal{F}\) is the smallest \(\sigma
	\)-algebra containing \(\mathcal{F}\times \mathcal{F}\)) whose
	marginals are \(\mathbb{P}\) and \(\mathbb{P}'\), i.e. \begin{equation}
		\mathbb{P}=\widehat{P}\circ \pi^{-1},\quad
		\mathbb{P}'=\widehat{P}\circ \pi'^{-1},
	\end{equation}
	where \(\pi \) is the left-projection and \(\pi'\) is the
	right-projection, defined by \begin{equation}
		\pi(x,x')=x,\quad \pi'(x,x')=x'
	\end{equation}
	for all \((x,x')\in \Omega\times \Omega\).
\end{defi}

We can define an analogous concept for random variables.

\begin{defi}[Coupling of random variables]
	A coupling of two random variables \(X\) and \(X'\) taking values in
	\((\Omega,\mathcal{F})\) is any pair of random variables
	\((\widehat{X},\widehat{X}')\) taking values in
	\((\Omega\times\Omega,\mathcal{F}\otimes{} \mathcal{F})\) whose
	marginals have the same distribution as \(X\) and \(X'\), i.e., \begin{equation}
		\widehat{X}\overset{D}{=}X,\quad \widehat{X}'\overset{D}{=}X',
	\end{equation}
	with \(\overset{D}{=}\) denoting equality in distribution.
\end{defi}

The following definition is inspired by the definition presented in
\cite[p.12]{coupling}.

\begin{defi}[Coupling time]
	Let \(x=(x(t))_{t\ge 0}\) and \(y=(y(t))_{t\ge 0}\) be two stochastic
	processes taking values in \((\Omega^{\mathbb{R}_{\ge
			0}},\mathcal{F}^{\otimes{} \mathbb{R}_{\ge 0}})\), let \((\widehat{x},\widehat{y})\) be a coupling of \(x\) and
	\(y\). Then we define \begin{equation}
		T=\inf\{t\ge 0:x(s)=y(s)\text{ for all }s\ge t\},
	\end{equation}
	which is the \textbf{coupling time} of \(x\) and \(y\).
\end{defi}

\subsection{Generalizing It\^o's formula and the time change argument to the local
	case}
Since we are interested in SDEs with blow-up solutions, the It\^o formula is
typically not immediately applicable. However, we can still apply the It\^o
formula to exploding SDEs. It then holds until an explosion time
\(S_\infty{}\). This is probably a known result which we could not
find in the literature, so we are presenting a proof here.

\begin{satz}[It\^o's formula for local solutions]
	\label{thm:itosformula_explosions}
	Let \(b\) and \(\sigma \) be adapted processes and let \(x\) be a solution of \begin{equation}
		\,\mathrm{d}x(t)=b(t)\,\mathrm{d}t+\sigma(t)\,\mathrm{d}W(t),
	\end{equation}
	up to an explosion time \(S_\infty{}\).
	Let \(f:\mathbb{R}\to\mathbb{R}\) be a \(C^2\)-function, then we have
	\(\mathbb{P}\)-a.s. \begin{align}
        \,\mathrm{d}f(x(t))=&\left[f'(x(t))b(t)+\frac{1}{2}f''(x(t))\sigma^2(t)\right]\,\mathrm{d}t\\
        +&f'(x(t))\sigma(t)\,\mathrm{d}W(t),
	\end{align}
	for all \(t\in[0,S_\infty{})\). It\^o's formula thus applies until
	the explosion time.
\end{satz}

\begin{bew}
	Consider the process satisfying
	\begin{align}
		\,\mathrm{d}x^n(t) & =b^n(t)\,\mathrm{d}t+\sigma^n(t)\,\mathrm{d}W(t),\intertext{where
			we define \(S_n:=\inf\{t>0:|x(t)|>n\}\),}
		b^n(t)             & = \begin{cases}
			                       b(t), & t<S_n,    \\
			                       0,    & t\ge S_n,
		                       \end{cases}                                               \\
		\sigma^n(t)        & = \begin{cases}
			                       \sigma(t), & t<S_n,    \\
			                       0,         & t\ge S_n.
		                       \end{cases}
	\end{align}
	It is clear that we have \(S_n\uparrow S_\infty{}\). Solutions \(x^n\) are continuous semimartingales, which means
    that It\^o's formula as stated in \cite[Theorem 3.3.3,
    p.149]{karatzas_shreve}
	applies. This yields that for all \(n\in\mathbb{N}\) we have:
	\begin{align}
        &\mathrm{d}f(x^n(t))=\left[f'(x^n(t))b^n(t)+\frac{1}{2}f''(x^n(t))(\sigma^n)^2(t)\right]\,\mathrm{d}t\\
        +&f'(x^n(t))\sigma^n(t)\,\mathrm{d}W(t).
	\end{align}
	As we have \(x^n=x,b^n=b\) and \(\sigma^n=\sigma\) for all
	\(t\in[0,S_n)\), we can conclude
	that we have in the limit\begin{align}
        &\mathrm{d}f(x(t))=\left[f'(x(t))b(t)+\frac{1}{2}f''(x(t))\sigma^2(t)\right]\,\mathrm{d}t\\
        +&f'(x(t))\sigma(t)\,\mathrm{d}W(t)
	\end{align}
	\(\mathbb{P}\)-a.s. for all \(t\in[0,S_\infty{})\).
	This concludes the argument.
\end{bew}

\begin{satz}[Time change argument for local solutions]
	\label{thm:time_change_fixed}
	Let \((\Omega,\mathcal{F},\mathbb{P}),(\mathcal{F}_t)_{t\ge 0}\) be a
	filtered probability space, and let \(S_{\text{stop}}\) be a stopping
	time.
	Let \((\sigma(t))_{t\in[0,S_{\mathrm{stop}})}\) be an adapted measurable
	process with \begin{equation}
		\mathbb{E} \int_0^{t\land{}S} \sigma(s)^2\,\mathrm{d}s<\infty{},\quad\text{for
			all }t\ge 0,
	\end{equation}
	for all stopping times \(0\le S<S_{\text{stop}}\) and with
	\(\sigma^2(t)>0\) for all \(t\in[0,S_{\mathrm{stop}})\) almost surely.
	Defining \begin{equation}
		T(t):=\int_0^t\sigma^2(s)\,\mathrm{d}s
	\end{equation}
	for \(t\in[0,S_{\mathrm{stop}})\), we get the following equality: \begin{equation}
		\int_0^t \sigma(s)\,\mathrm{d}W(s)=B(T(t)),\quad\text{a.s. for all
		}t\in[0,S_{\mathrm{stop}}),
	\end{equation}
	where \(B\) is a Brownian motion with respect to a suitable filtration
	on \((\Omega,\mathcal{F},\mathbb{P})\)..
\end{satz}

\begin{bew}
	Define \begin{equation}
		S_n:=\inf\left\{t\in(0,S_{\text{stop}}):\sigma(s)\notin\left(\frac{1}{n},n\right)\right\}.
	\end{equation}
	Then we define \begin{equation}
		\sigma^n(s)= \begin{cases}
			\sigma(s),   & s< S_n,   \\
			\sigma(S_n), & s\ge S_n.
		\end{cases}
	\end{equation}
	It is clear that \(\sigma^n\) fulfills the conditions of the time
    change argument \cite[Theorem IV.34.1]{rogers_williams}.
    We thus have \begin{equation}
		\int_0^t \sigma^n(s)\,\mathrm{d}W(s)=B_n(T(t))
	\end{equation}
	almost surely for all \(t\ge 0\).
	It is clear that up until \(S_n \land{}S_m\) the Brownian motions
	\(B_n\) and \(B_m\) equal each other. Since we have this uniqueness,
	we simply denote them by \(B\). This means that we have \begin{equation}
		\int_0^{t\land{}S_n}\sigma(s)\,\mathrm{d}W(s)=B(T(t\land{}S_n))
	\end{equation}
	for all \(t>0\), \(n\in\mathbb{N}\). Note that
	\(S_{\text{stop}}=\lim_{n\to \infty{}}S_n\), thus
	the following holds for all \(t\in[0,S_{\mathrm{stop}})\):
	\begin{equation}
		\int_0^t \sigma(s)\,\mathrm{d}W(s)=B(T(t)).
	\end{equation}
\end{bew}

%% file: 40_delay.tex
\section{Theory}
\label{sect:delayed_models}
This allows for analyzing explosion and extinction behavior of SDDEs
by lifting the results from \cite[Section 4]{masters}.
\subsection{Unperturbed case: Known results for DDEs}
\label{sect:delayed_models_known_results}

We begin our investigations by presenting the findings of Ishiwata and
Nakata in \cite{jap}.
Let \(F\) be a continuous function from \(\mathbb{R}^2\) to
\(\mathbb{R}\). We consider the DDE \begin{equation}
x'(t)=F(x(t),x(t-1))
\label{eq:jap:dde}
\end{equation}
with the initial condition \begin{equation}
x(t)=\phi(t),\ t\in[-1,0],
\end{equation}
where \(\phi:[-1,0]\to \mathbb{R}\) is continuous. Ishiwata and Nakata prove
that under certain regularity conditions imposed on \(F\), the
explosion behavior of \eqref{eq:jap:dde} is closely related to the
explosion behavior of the ODE \begin{equation}
x'(t)=F(x(t),a),
\label{eq:jap:ode}
\end{equation}
where \(a\in\mathbb{R}\) is some parameter. Even without requiring any
special regularity properties, the converse holds true \cite[Theorem 3.1]{jap}.

\begin{satz}[Relating ODE blow-up to DDE blow-up]
\label{thm:jap1}
If there exists an \(a\in\mathbb{R}\) such that the ODE
\eqref{eq:jap:ode} has a blow-up solution, then the DDE
\eqref{eq:jap:dde} also has blow-up solutions.
\end{satz}

This theorem is proven easily by defining a suitable initial function \(\phi
\).
The converse only holds when a further
condition is imposed. If
for any closed bounded interval
\(I\subset{}\mathbb{R}\) there exist \(a_1,a_2\in\mathbb{R}\) such
that
\begin{equation}
(x,y)\in\mathbb{R}\times I\Rightarrow{}F(x,a_1)\le F(x,y)\le F(x,a_2),
\label{jap:cond}
\end{equation}
then we say that \(F\) fulfills \eqref{jap:cond}.
We continue by stating the following theorem \cite[Theorem
3.2]{jap}.\begin{satz}[Relating DDE blow-up to ODE blow-up]
\label{thm:jap2}
Suppose that \eqref{jap:cond} holds. If the DDE \eqref{eq:jap:dde} has a
blow-up solution, then there exists an \(a\in\mathbb{R}\) such that
the ODE \eqref{eq:jap:ode} also has blow-up solutions.
\end{satz}
\subsection{Generalizing to the perturbed case}
In this section we study SDDEs of the form \begin{align}
\label{eq:jap:sdde}
\,\mathrm{d}x(t)&=F(x(t),x(t-1))\,\mathrm{d}t+g(x(t))\,\mathrm{d}W(t),\
t>0,\\
x(t)&=\phi(t),\ t\in[-1,0].
\end{align}
where \(\phi:[-1,0]\to \mathbb{R}\) is continuous and \(g\) satisfies
\eqref{eq:K}. This is
assumed to hold from now on. Furthermore, in the interest of
simplifying existence and uniqueness considerations, we consider \(F\)
and \(g\) to be locally Lipschitz.
We want to expand the results of Section
\ref{sect:delayed_models_known_results} to the stochastically perturbed
case.
As before, we can relate the explosion behavior of this SDDE to the
explosion behavior of the SDE \begin{equation}
\,\mathrm{d}x(t)=F(x(t),a)\,\mathrm{d}t+g(x(t))\,\mathrm{d}W(t),
\label{eq:jap:sde}
\end{equation}
with \(x(0)=x_0\).
We begin
our exploration by generalizing Theorem \ref{thm:jap1} to the
stochastic case. 
\begin{satz}[Relating SDE blow-up to SDDE blow-up]
\label{thm:trivial_impl}
If there exists an \(a\in\mathbb{R}\) such that the corresponding solution of the SDE
\eqref{eq:jap:sde} blows up to \(\pm\infty{}\) with nonzero probability, then there
exists an initial function such that the solution of the SDDE
\eqref{eq:jap:sdde} blows up to \(\pm\infty{}\) before
\(t=1\) with nonzero probability.
\end{satz}

\begin{bew}
We shall provide an appropriate initial
function. Remember that \(F\) and \(g\) are locally Lipschitz.
Let \(x\) denote the solution of the SDE \eqref{eq:jap:sde}. This
solution \(x\) blows
up to \(\pm\infty{}\) with nonzero probability. Let \(T\) denote the explosion time of
this solution and let \(x_0=x(0)\) be the initial
condition. It is immediate that there exists an
interval \([n,n+1)\) with \(n\in \mathbb{Z}_{\ge 0}\) such that \(\mathbb{P}[T\in[n,n+1)]>0\) (see also
Lemma \ref{lem:sdde_to_finite_ints}) Let \([n,n+1) \) denote the
smallest interval in which this holds true.

If \(n=0\), there is some
\(\epsilon > 0\) such that \(\mathbb{P}[T\in [0,1-\epsilon)]>0\). Then
we define \begin{equation}
\psi(t)=
\begin{cases}
a, & t\in [-1,-\epsilon),\\
\frac{a-x_0}{-\epsilon }t+x_0,&t\in[-\epsilon,0].
\end{cases}
\end{equation}
Since the solution \(\overline{x}\) of the SDDE \eqref{eq:jap:sdde}
now has the same distribution as the solution \(x\) of the SDE
\eqref{eq:jap:sde} up until \(t=1-\epsilon\),
we see that \(\overline{x}\) also has a nonzero probability of
exploding to \(\pm\infty{}\) in finite time in that case.

If \(n\ne 0\), define
\begin{equation}
\psi(t)=x(t+n),\ t\in[-1,0].
\end{equation}
This is well-defined because the explosion probability is zero before
\(t=n\).
Let now \(\overline{x}\) denote the solution of the SDDE \eqref{eq:jap:sdde} with initial
function \(\psi \) instead of $\phi$. Clearly the distribution of \(\overline{x}\) in
\([0,1]\) is identical to the distribution of \(x\) in
\([n,n+1\land{}T]\), where \(T\) denotes the explosion time. This means that \(x\) has a nonzero probability of exploding to
\(\pm\infty{}\) before
\(t=1\), which concludes the proof.
\end{bew}

Before we proceed to the main result, we present some
observations.

\begin{lem}
\label{lem:sdde_to_finite_ints}
Given an SDDE \begin{equation}
\,\mathrm{d}x(t)=F(x(t),x(t-1))\,\mathrm{d}t+g(x(t))\,\mathrm{d}W(t),
\end{equation}
with initial condition \(x(t)=\phi(t)\) for \(t\in[-1,0]\) where \(g\)
fulfills \eqref{eq:K} and the coefficients $F$ and $g$ are locally Lipschitz, there
is a nonzero probability that a solution \(x\) blows up to \(\pm\infty{}\) in finite time if and only if there
exists some interval \([n,n+1)\) such that
\(\mathbb{P}(T\in[n,n+1))>0\).
\end{lem}

\begin{bew}
Assume the negation for the purpose of contradiction. As
\(\mathbb{P}(T\in[n,n+1))=0\) for all \(n\in \mathbb{N}\), we have \begin{align}
    0<\mathbb{P}(T\in[0,\infty{}))=&\mathbb{P}\left(T\in\bigcup_{n\ge
0}[n,n+1)\right)\\
            =&\sum_{n\ge
0}\mathbb{P}(T\in[n,n+1))=0,
\end{align}
which yields a contradiction.
The converse direction is trivial.
This proves the statement.
\end{bew}

\begin{lem}
\label{lem:sdde_to_finite_2}
Consider the SDDE \begin{align}
&\mathrm{d}x(t)=F(x(t),x(t-1))\,\mathrm{d}t+g(x(t))\,\mathrm{d}W(t),\\
&x(t)=\phi(t),\ t\in[-1,0],
\end{align}
where \(g\) fulfills \eqref{eq:K} and both coefficients $F$ and $g$ are locally
Lipschitz.
If this
SDDE has a nonzero probability of exploding to \(\pm\infty{}\) in finite time, then there
exists an initial stochastic continuous function \(\psi\) such that the SDDE \begin{align}
&\mathrm{d}x(t)=F(x(t),x(t-1))\,\mathrm{d}t+g(x(t))\,\mathrm{d}W(t),\\
&x(t)=\psi(t),\ t\in[-1,0]
\end{align}
has a nonzero probability of exploding to \(\pm\infty{}\) before \(t=1\).
\end{lem}

\begin{bew}
The proof relies on Lemma \ref{lem:sdde_to_finite_ints}. That
lemma provides an interval \([n,n+1)\) with \(n\in\mathbb{Z}_{\ge 0}\) such that \(x\) has a nonzero
probability of exploding to \(\pm\infty{}\) in \([n,n+1)\). Let
\([n,n+1)\) denote such and interval nearest to zero. If \(n=0\), we are
immediately done by picking \(\psi=\phi\). Assume thus that \(n>0\). Let now $\phi(t)$ be a $C[-1,0]$-valued random variable with the same distribution as $t\mapsto x(n+t)$ and independent of $W$.
This is well-defined, because solutions do not explode in finite time
before \(t=n\) with probability one.
Because we presume solutions of \(x\) to be unique up to the explosion
time, it is clear that the SDDE
\begin{align}
&\mathrm{d}x(t)=F(x(t),x(t-1))\,\mathrm{d}t+g(x(t))\,\mathrm{d}W(t),\\
&x(t)=\psi(t),\ t\in[-1,0]
\end{align}
has a nonzero probability of exploding to \(\pm\infty{}\) before \(t=1\).
This concludes the argument.
\end{bew}

\begin{lem}
\label{lem:sdde_reduction}
Consider an SDDE of the form \begin{equation}
\,\mathrm{d}x(t)=F(x(t),x(t-1))\,\mathrm{d}t+g(x(t))\,\mathrm{d}W(t),\
t\ge 0,
\end{equation}
with initial function \(x(t)=\phi(t),\ t\in[-1,0]\) with locally
Lipschitz coefficients $F$ and $g$.
If \(F\) fulfills
\eqref{jap:cond} and \(g\) fulfills the condition \eqref{eq:K} from Theorem
\ref{thm:watanabe_comparison_sde}, then there exist constants
\(a_1,a_2\) such that the solutions of the SDEs \begin{align}
    \mathrm{d}x_1(t)&=F(x_1(t),a_1)\,\mathrm{d}t+g(x_1(t))\,\mathrm{d}W(t),\\
    x_1(0)&=\phi(0),\\
    \mathrm{d}x_2(t)&=F(x_1(t),a_2)\,\mathrm{d}t+g(x_2(t))\,\mathrm{d}W(t),\\
    x_2(0)&=\phi(0),
\end{align}
fulfill \begin{equation}
x_1(t)\le x(t)\le x_2(t)
\end{equation}
up to \(T\land{}1\)
with probability one, where \(T\) denotes the explosion time.
\end{lem}

\begin{bew}
We proceed similarly as in the proof of \cite[Theorem 3.2]{jap}. Since
\(\phi \) is continuous and defined on a compact interval, the
condition \eqref{jap:cond} provides \(a_1,a_2\) such that \begin{equation}
F(z,a_1)\le F(z,\phi(t-1))\le F(z,a_2)
\end{equation}
for all \(t\in[0,1]\cap [0,T]\), $z\in\mathbb R$. 
With these $a_1$ and $a_2$ we associate the SDEs
\begin{align}
\,\mathrm{d}x_1(t)=F(x_1(t),a_1)\,\mathrm{d}t+g(x_1(t))\,\mathrm{d}W(t),\\
\,\mathrm{d}x_2(t)=F(x_2(t),a_2)\,\mathrm{d}t+g(x_2(t))\,\mathrm{d}W(t).
\end{align}
Since the solution is unique, \(g\) fulfills \eqref{eq:K} and the
starting values are all equal,
Theorem
\ref{thm:watanabe_comparison_sde}
then states that we have \begin{equation}
x_1(t)\le x(t)\le x_2(t)
\end{equation}
for all \(t\in[0,1]\cap[0,T]\) almost surely. This concludes the
proof of the lemma.
\end{bew}

\begin{rem}
\label{rem:stoch_init}
Note that the statement of Lemma \ref{lem:sdde_reduction} also holds for stochastic initial conditions
\(\phi \) which are independent of $W$. If the equation is well-defined, then one can factor the
probability space
into a direct product of the space where the initial condition
\(\phi\) lives one and into one that carries $W(t)$ for \(t> 0\).
\end{rem}

We now proceed to the main result of this section.

\begin{satz}[Reducing explosion behavior of SDDE to SDE]~\\
\label{thm:sdde_explosion_equivalence}
Consider an SDDE of the form 
\begin{align}
&\mathrm{d}x(t)=F(x(t),x(t-1))\,\mathrm{d}t+g(x(t))\,\mathrm{d}W(t),\
t\ge 0,\\
&x(t)=\phi(t),\ t\in[-1,0],
\end{align}
where \(F\) and \(g\) are locally Lipschitz, \(F \) fulfills
\eqref{jap:cond} and \(g\) fulfills \eqref{eq:K}.
There exists an initial function \(\phi \) such that solutions \(x\) have a nonzero probability of exploding to \(\pm\infty{}\) in finite
time if and only if there exists an \(a\in\mathbb{R}\) and an initial
condition \(x_0\) such that the
SDE \begin{align}
&\mathrm{d}x(t)=F(x(t),a)\,\mathrm{d}t+g(x(t))\,\mathrm{d}W(t),\\
&x(0)=x_0
\end{align}
has a nonzero probability of exploding to \(\pm\infty{}\) before \(t=1\). The statement
also holds if we replace ``exploding before \(t=1\)'' by ``exploding
in finite time''.
\end{satz}

\begin{bew}
The implication ``\(\Leftarrow\)'' was covered in Theorem
\ref{thm:trivial_impl}, since if the solution of the SDE blows up to \(\pm\infty{}\)
before \(1\) with nonzero probability, it also blows up to \(\pm\infty{}\) in finite time
with nonzero probability. Assume thus that the
SDDE has a nonzero probability of exploding in finite time. Lemma
\ref{lem:sdde_to_finite_2} gives a new SDDE 
\begin{align}
&\mathrm{d}x(t)=F(x(t),x(t-1))\,\mathrm{d}t+g(x(t))\,\mathrm{d}W(t),\
t\ge 0,\\
&x(t)=\psi(t),\ t\in[-1,0],
\end{align}
which has a nonzero probability of exploding to \(\pm\infty{}\) before \(t=1\). Applying
Lemma \ref{lem:sdde_reduction} (see also Remark \ref{rem:stoch_init}) now yields constants \(a_1,a_2\) such
that the solutions
\begin{align}
\,\mathrm{d}x_1(t)=F(x_1(t),a_1)\,\mathrm{d}t+g(x_1(t))\,\mathrm{d}W(t),\\
\,\mathrm{d}x_2(t)=F(x_1(t),a_2)\,\mathrm{d}t+g(x_2(t))\,\mathrm{d}W(t),
\end{align}
fulfill \begin{equation}
x_1(t)\le x(t)\le x_2(t)
\end{equation}
with probability one up to the explosion time. If explosion of $x$ occurs to
\(+\infty{}\), then \(x_2\) also explodes. We argue similarly for the case
of explosion to \(-\infty{}\). This concludes the proof.
\end{bew}

\subsection{Example}
\begin{beisp}[Delay prevents extinction]
\label{ex:delay}
Consider the SDDE population model of the form \begin{align}
\label{eq:sdde:simple_delay_pop_model}
\,\mathrm{d}y(t)&=[ay(t-1)-by(t)]\,\mathrm{d}t+y(t)^p\,\mathrm{d}W(t),\ t\ge 0,\\
y(t)&=\phi(t)>0,\ t\in[-1,0],
\end{align}
with \(p\in(0,1)\).
It is known that we have unique solutions up to a stopping time \(T\), which
in this case is \(\lim_{n\to \infty{}}\inf\{t>0:t\notin(\frac{1}{n},n)\}\).
We begin by examining its symbolic transformation with It\^o's
formula as formulated in Theorem \ref{thm:itosformula_explosions},
using \(f\equiv -\log\). In the usual manner, this will be made rigorous by
transforming back later. We thus consider the SDDE
\begin{align}
    \,\mathrm{d}x(t)=\Big[&-e^{x(t)}\left(ae^{-x(t-1)}-be^{-x(t)}\right)\\
                           &+\frac{1}{2}e^{2(1-p)x(t)}\Big]\,\mathrm{d}t+e^{(1-p)x(t)}\,\mathrm{d}W(t).
\end{align}
Note that the coefficients are locally Lipschitz, and we thus have
unique solutions up to an explosion time. We can now apply Theorem
\ref{thm:sdde_explosion_equivalence}. This means that it is sufficient
to consider the explosion behavior of the SDE \begin{align}
    \,\mathrm{d}x(t)=\Big[&-e^{x(t)}\left(ae^{-C}-be^{-x(t)}\right)\\
&+\frac{1}{2}e^{2(1-p)x(t)}\Big]\,\mathrm{d}t+e^{(1-p)x(t)}\,\mathrm{d}W(t),
\end{align}
where \(C\in\mathbb{R}\) is some constant. We can now rigorously
transform this equation back onto \(\mathbb{R}_+\), using Theorem
\ref{thm:itosformula_explosions} with \(f\equiv e^{-\cdot}\). This
yields the SDE \begin{align}
    \,\mathrm{d}y(t)&=\left[ae^{-C}-by(t)\right]\,\mathrm{d}t+y(t)^p\,\mathrm{d}W(t),\\
    y(0)&=y_0>0,
\end{align}
up until the stopping time \(T\). If \(p<\frac{1}{2} \), it follows
immediately from Theorem 
\cite[Theorem 4.21, p.54]{masters} that
we have a nonzero probability of extinction. If \(p>\frac{1}{2}\),
then it follows from \cite[Theorem 4.28, p.57]{masters} 
that solutions
never hit zero in finite time. In this way,
extinction results for SDEs carry over to the delayed case.
\end{beisp}

%% file: 60_summary.tex
\section{Conclusions}
After presenting the results of Ishiwata and Nakata \cite{jap}, we
generalized their findings to the stochastically perturbed case. In
particular, we showed that explosion 
properties of SDDEs whose
coefficients fulfill
\eqref{jap:cond} and \eqref{eq:K} can be reduced to explosion
properties of corresponding SDEs by replacing the delay term
\(x(t-\tau)\) by a constant \(a\in\mathbb{R}\). To the knowledge of
the author
this is a novel result which allows
extensive description of explosion and extinction dynamics of SDDEs.
Finally, we employed this technique to analyze
extinction behavior of the simple population model
\(\mathrm{d}y(t)=[ay(t-1)-by(t)]\,\mathrm{d}t+y(t)^p\,\mathrm{d}W(t)\),
proving that for \(a,b,y(0)>0\) extinction can occur for
\(p<\frac{1}{2}\) but not for \(p>\frac{1}{2}\).
We generalized the result by \cite{jap} to stochastically perturbed equations,
using a comparison theorem by Ikeda and Watanabe \cite{watanabe}. This, combined with \cite[Section 4]{masters}, allows analyzing
extinction properties of SDDEs.
We demonstrated that the extinction dynamics of sufficiently
regular delayed models can be effectively reduced to those of
non-delayed models, replacing the delay term with a constant value.

\section*{Acknowledgements}
I wish to thank Dr. Onno van Gaans for supervising me during the preparation of
this paper.